\documentclass[12pt,a4paper,leqno]{article}
\usepackage{amsmath,amssymb,amsthm,bm,latexsym}

\def\R{\mathbb{R}}
\def\Q{\mathbb{Q}}

\def\N{\mathcal{N}}
\def\M{\mathcal{M}}
\def\K{\mathcal{K}}

\def\add#1{{\operatorname{\textup{\textsf{add}}}}(#1)}
\def\non#1{{\operatorname{\textup{\textsf{non}}}}(#1)}
\def\cov#1{{\operatorname{\textup{\textsf{cov}}}}(#1)}
\def\cof#1{{\operatorname{\textup{\textsf{cof}}}}(#1)}
\def\fb{\mathfrak{b}}
\def\fc{\mathfrak{c}}
\def\fd{\mathfrak{d}}

\def\fu{\mathfrak{u}}
\def\ssm{\smallsetminus}

\def\Pow#1{\mathcal{P}(#1)}
\def\restrictedto%
{\mathchoice%
    {\!\upharpoonright\!}
    {\!\upharpoonright\!}
    {\upharpoonright}
    {\upharpoonright}
}

\def\cf{\operatorname{cf}}
\def\op#1#2{\langle #1,#2\rangle}

\def\st{\mathchoice{:}{:}{\,:\,}{\,:\,}}

\def\size#1{\lvert#1\rvert}

\def\domn{\leq^*}
\def\aincl{\subseteq^*}

\newcommand{\splitnodes}[2][{}]%
	{\operatorname{\textup{\textsf{split}}}_{#1}({#2})}

\def\abs#1{\left\lvert{#1}\right\rvert}
\def\V{\mathbf{V}}

\def\cR{\mathcal{R}}
\def\cl{\operatorname{cl}}	
\def\Ball{\operatorname{B}}
\def\diam{\operatorname{diam}}
\def\higson#1#2{\overline{#1}^{#2}}
\def\smirnov#1#2{u_{#2}{#1}}
\def\stonecech#1{\beta{#1}}
\def\PM{\operatorname{PM}}
\def\Metric{\operatorname{M}}
\def\cptequiv{\simeq}

\def\fl{\mathfrak{l}}
\def\fsa{\mathfrak{sa}}
\def\fsp{\mathfrak{sp}}
\def\fst{\mathfrak{st}}
\def\fha{\mathfrak{ha}}
\def\fhp{\mathfrak{hp}}
\def\fht{\mathfrak{ht}}
\def\fpp{\mathfrak{pp}}

\def\cptsep#1#2#3{{#1}\parallel{#2}\,\,\,(#3)}
\def\cptnonsep#1#2#3{{#1}\not\,\parallel{#2}\,\,\,(#3)}

\def\SIGMA{\bm{\Sigma}}
\def\PI{\bm{\Pi}}

\theoremstyle{plain}
\newtheorem{thm}{Theorem}[section]
\newtheorem{lem}[thm]{Lemma}
\newtheorem{cor}[thm]{Corollary}
\newtheorem{prop}[thm]{Proposition}

\newtheorem{claim}{Claim}

\newtheorem{q}[thm]{Question}
\newtheorem{ex}[thm]{Example}
\theoremstyle{definition}
\newtheorem{defn}[thm]{Definition}
\theoremstyle{remark}

\newtheorem{remark}{Remark}

\begin{document}

\title{How many miles to $\beta\omega$? 
	--- Approximating 
	$\beta\omega$ by 
	metric-dependent compactifications
	}
	
\author{Masaru Kada%
    \thanks{%
	Supported by 
	Grant-in-Aid for Young Scientists (B) 14740058, MEXT.} 
	\and
	Kazuo Tomoyasu%
	\thanks{%
	Supported by 
	Grant-in-Aid for Young Scientists (B) 14740057, MEXT.} 
	\and 
	Yasuo Yoshinobu
	\thanks{%
	Supported by
	Grant-in-Aid for Scientific Research (C)(2) 15540115, JSPS.}
	}
\date{}
\maketitle

\begin{abstract}
It is known that 
the Stone-\v{C}ech compactification $\stonecech{X}$ 
of a non-compact 
metrizable space $X$ 
is approximated by the collection of Smirnov compactifications of $X$ 
for all compatible metrics on $X$. 
We investigate the smallest cardinality 
of a set $D$ of compatible metrics 
on the countable discrete space $\omega$ such that, 
$\stonecech{\omega}$ is approximated by Smirnov compactifications 
for all metrics in $D$, 
but any finite subset of $D$ does not suffice. 
We also study the corresponding cardinality 
for Higson compactifications.


\par
\vspace{12pt}
\par
\vspace{12pt}
\end{abstract}

\section{Notations and definitions}\label{sec:intro}

\subsection{Compactifications of topological spaces}\label{subsec:compact}

Let $X$ be a completely regular Hausdorff topological space. 
A \emph{compactification $\alpha X$ of $X$} 
is a compact Hausdorff space 
which contains $X$ as a dense subspace. 
For compactifications $\alpha X$ and $\gamma X$ 
of a non-compact space $X$, 
we write $\alpha X\leq \gamma X$ 
if there is a continuous surjection $f\st \gamma X\to\alpha X$ 
such that 
$f\restrictedto X$ is the identity map on $X$.
If such an $f$ can be chosen to be a homeomorphism, 
we say $\alpha X$ and $\gamma X$ are \emph{equivalent} 
and denote this by writing $\alpha X\cptequiv\gamma X$. 

Let $\K(X)$ denote the class of compactifications of $X$. 
When 
we identify equivalent compactifications 
and regard $\K(X)$ as the collection of equivalence classes, 
we may regard $\K(X)$ as a set, 
and then the order structure $(\K(X),{\leq})$ 
is a complete upper semilattice 
whose largest element is the Stone-\v{C}ech compactification $\stonecech{X}$.

Let $C^*(X)$ denote the set of all bounded continuous functions 
from $X$ to $\R$. 
$C^*(X)$ is a topological ring 
with respect to pointwise addition, pointwise multiplication, 
and the uniform norm topology (see \cite[2M]{GJ:rings}).
A subring $R$ of $C^*(X)$ is called \emph{regular} if $R$ is closed,
contains all constant functions, and generates the topology of $X$.
Let $\cR(X)$ denote the class of regular subrings of $C^*(X)$.
Then it is known that $(\K(X), {\leq})$ is isomorphic to
$(\cR(X), {\subseteq})$, by mapping each $\alpha X\in\K(X)$
to the set of bounded continuous functions from $X$ to $\R$ 
which are continuously extended over $\alpha X$, which we will
denote as $C_{\alpha X}$ (cf.\ \cite[Theorem 3.7]{BaY:ring}, \cite[Theorem 2.5]{Chan:cpt}).
In particular, the Stone-\v{C}ech compactification $\beta X$
corresponds to the whole $C^*(X)$.
(See \cite{Chan:cpt,GJ:rings} for more details.)

We introduce the following notation. 
For a compactification $\alpha X$ of $X$ 
and two closed subsets $A,B$ of $X$, 
we write $\cptsep{A}{B}{\alpha X}$ 
if $\cl_{\alpha X}A\cap\cl_{\alpha X}B=\emptyset$, 
and otherwise we write 
$\cptnonsep{A}{B}{\alpha X}$. 

The following lemma 
is well-known.

\begin{lem}\label{lem:charcptsep}
For a compactification $\alpha X$ of a space $X$ 
and closed subsets $A,B$ of $X$, 
the following are equivalent\/\textup{:}
\begin{enumerate}
\item $\cptsep{A}{B}{\alpha X}$. 
\item There is an $f\in C_{\alpha X}$ 
	such that 
	$f''A=\{1\}$ and $f''B=\{0\}$. 
\item There are $g\in C_{\alpha X}$ 
	and $a,b\in\R$ with $a>b$ 
	such that 
	$g(x)\geq a$ for all $x\in A$ and 
	$g(x)\leq b$ for all $x\in B$. 
\end{enumerate}
\end{lem}

Note that, for compactifications $\alpha X$, $\gamma X$ of a space $X$,
$\alpha X\leq\gamma X$ holds if and only if, 
$\cptsep{A}{B}{\alpha X}$ implies $\cptsep{A}{B}{\gamma X}$
for any $A,B\subseteq X$.
In particular,
$\alpha X\cptequiv\stonecech{X}$
if and only if $\cptsep{A}{B}{\alpha X}$ holds
for any pair of disjoint closed subsets $A,B$ of $X$.

%

The following lemma 
plays an important role 
in the following sections. 

\begin{lem}\label{lem:cptsepcompact}
Suppose that $\mathcal{C}$ is a set of compactifications of a space $X$. 
For closed sets $A,B$ of $X$, 
the following are equivalent 
\textup{(}$\sup$ is taken in the lattice $(\mathcal{K}(X),{\leq})$\textup{):} 
\begin{enumerate}
\item $\cptsep{A}{B}{\sup\mathcal{C}}$. 
\item $\cptsep{A}{B}{\sup\mathcal{F}}$ 
	for some nonempty finite subset $\mathcal{F}$ of\/ $\mathcal{C}$. 
\end{enumerate}
\end{lem}

\begin{proof}
Let $R=\langle\bigcup\{C_{\alpha X}\st\alpha X\in\mathcal{C}\}\rangle$ 
(where $\langle Z\rangle$ denotes the subring of $C^*(X)$ 
generated by $Z\subseteq C^*(X)$) 
and $C=\cl R$. 
Then we have $C=C_{\sup\mathcal{C}}$. 

Suppose that $A,B\subseteq X$ are closed 
and $\cptsep{A}{B}{\sup\mathcal{C}}$. 
By Lemma \ref{lem:charcptsep}, 
there is an $f\in C$ such that $f''A=\{1\}$ and $f''B=\{0\}$. 
Since $C=\cl R$, 
there is a sequence of functions in $R$ 
which uniformly converges to $f$, 
and so there is a $g\in R$ such that 
$g(x)\geq\frac{2}{3}$ if $x\in A$ and $g(x)\leq\frac{1}{3}$ if $x\in B$. 
By the definition of $R$, 
there is a finite set $\mathcal{F}\subseteq\mathcal{C}$ 
such that 
$g$ is obtained by additions and multiplications 
of functions from 
$\bigcup\{C_{\alpha X}\st\alpha X\in\mathcal{F}\}$, 
and hence 
$g\in R'=\langle\bigcup\{C_{\alpha X}\st\alpha X\in\mathcal{F}\}\rangle$. 
Then 
$R'\subseteq\cl R'=C_{\sup\mathcal{F}}$, 
and by Lemma \ref{lem:charcptsep}, 
we have $\cptsep{A}{B}{\sup\mathcal{F}}$. 
\end{proof}


\subsection{Cardinal invariants of the reals}\label{subsec:ci}

In this subsection, 
we review definitions and known results about 
cardinal invariants of the reals, 
which appear in the following sections. 

We use standard notations and basic facts about set theory, 
including descriptive set theory. 
We refer the readers to 
\cite{BaJ:set}, \cite{Je:set}, \cite{Ke:descriptive} or \cite{Ku:set}  
for undefined set-theoretic notions. 
For $X,Y\subseteq\omega$, 
we write $X\aincl Y$ 
if $X\ssm Y$ is finite. 
For $f,g\in\omega^\omega$, 
$f\domn g$ 
if for all but finitely many $n<\omega$ we have $f(n)\leq g(n)$.

We will often use the expression 
``$\kappa$ is the smallest cardinality of a set $X$ which satisfies ...,'' 
without assuming the existence of such an $X$. 
Here we make the following notational convention: 
If there is no such $X$, 
we write $\kappa=\infty$ 
and regard $\lambda<\infty$ for any cardinal $\lambda$. 

Let $2^\omega$ be the Cantor space equipped with the usual product measure. 
Let $\M$ denote the collection of meager subsets of $2^\omega$, 
and $\N$ the collection of measure zero subsets of $2^\omega$. 
Both $\M$ and $\N$ are countably complete ideals on $2^\omega$ 
and contain all singletons. 

For a collection $\mathcal{X}$ of subsets of $2^\omega$, 
we define the following four cardinal coefficients. 
\begin{enumerate}
\item $\cov{\mathcal{X}}$
	is the smallest cardinality of a set $\mathcal{X}'\subseteq\mathcal{X}$ 
	such that $\bigcup\mathcal{X}'=2^\omega$. 
\item $\non{\mathcal{X}}$
	is the smallest cardinality of a set $Z\subseteq 2^\omega$ 
	such that $Z\not\subseteq X$ holds for every $X\in\mathcal{X}$. 
\item $\add{\mathcal{X}}$
	is the smallest cardinality of a set $\mathcal{X}'\subseteq\mathcal{X}$ 
	such that 
	$\bigcup\mathcal{X}'\not\subseteq X$ holds for every $X\in\mathcal{X}$. 
\item $\cof{\mathcal{X}}$
	is the smallest cardinality of a set $\mathcal{X}'\subseteq\mathcal{X}$ 
	such that 
	for every $Y\in\mathcal{X}$ there is an $X\in\mathcal{X}'$ 
	satisfying $Y\subseteq X$. 
\end{enumerate}

Clearly, 
$\mathcal{X}\subseteq\mathcal{Y}$ 
implies 
$\cov{\mathcal{X}}\geq\cov{\mathcal{Y}}$ 
and 
$\non{\mathcal{X}}\leq\non{\mathcal{Y}}$. 

$\fb$ is the smallest cardinality 
of an unbounded subset of a partially ordered set $(\omega^\omega,\domn)$. 
$\fd$ is the smallest cardinality 
of a cofinal subset of $(\omega^\omega,\domn)$. 

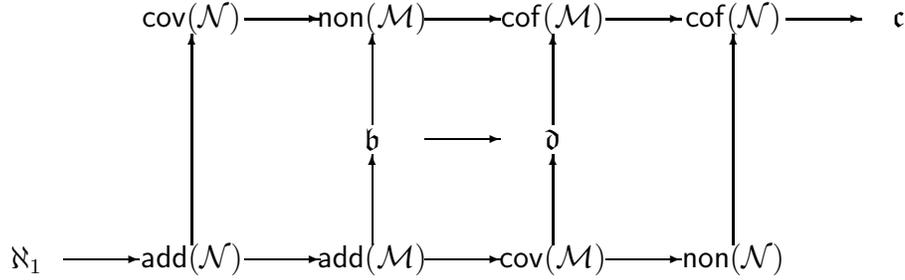
\begin{figure}
\begin{center}
\setlength{\unitlength}{0.2000mm}
\begin{picture}(636.0000,180.0000)(0,10)
\thinlines
\put(535,180){\vector(1,0){50}}
\put(415,180){\vector(1,0){50}}
\put(175,180){\vector(1,0){50}}
\put(295,180){\vector(1,0){50}}
\put(295,100){\vector(1,0){50}}
\put(415,20){\vector(1,0){50}}
\put(295,20){\vector(1,0){50}}
\put(175,20){\vector(1,0){50}}
\put(55,20){\vector(1,0){50}}
\put(140,30){\vector(0,1){140}}
\put(500,30){\vector(0,1){140}}
\put(380,30){\vector(0,1){60}}
\put(380,110){\vector(0,1){60}}
\put(260,110){\vector(0,1){60}}
\put(260,30){\vector(0,1){60}}
\put(470,170){\makebox(60,20){$\cof\N$}}
\put(470,10){\makebox(60,20){$\non\N$}}
\put(350,170){\makebox(60,20){$\cof\M$}}
\put(350,90){\makebox(60,20){$\fd$}}
\put(350,10){\makebox(60,20){$\cov\M$}}
\put(230,10){\makebox(60,20){$\add\M$}}
\put(230,90){\makebox(60,20){$\fb$}}
\put(230,170){\makebox(60,20){$\non\M$}}
\put(110,170){\makebox(60,20){$\cov\N$}}
\put(110,10){\makebox(60,20){$\add\N$}}
\put(10,10){\makebox(40,20){$\aleph_1$}}
\put(590,170){\makebox(40,20){$\fc$}}
\end{picture}
\end{center}
\caption{Cicho\'n's diagram}
\label{fig:cd}
\end{figure}

Figure~\ref{fig:cd}, 
which is known as ``Cicho\'n's diagram'' 
(\cite{Fr:cichon}), 
illustrates 
the relations among $\fb$, $\fd$ 
and cardinals defined from $\M$ and $\N$, 
where $\kappa\rightarrow\lambda$ in the diagram 
means that 
the inequality $\kappa\leq\lambda$ 
is provable in ZFC\@. 
Martin's axiom implies 
that 
all cardinals from the diagram, except for $\aleph_1$, equal to $\fc$. 
It is known that 
Cicho\'n's diagram is complete in the sense that 
every provable inequality is given by a chain of arrows in the diagram. 
(See \cite{BaJ:set} for details.)

Throughout this paper, 
an \emph{ultrafilter} means a non-principal ultrafilter on $\omega$. 
For an ultrafilter $\mathcal{U}$ and $\mathcal{F}\subseteq\mathcal{U}$, 
we say \emph{$\mathcal{F}$ generates $\mathcal{U}$}\ 
if for every $X\in\mathcal{U}$ 
there is a $Y\in\mathcal{F}$ such that $Y\subseteq X$. 
$\fu$ is the smallest cardinality 
of a subset of $\Pow{\omega}$ which generates an ultrafilter. 
Clearly we have $\fu\leq\fc$. 
%
%
It is known that $\cov\M\leq\fu$, $\cov\N\leq\fu$ 
(see \cite{Vo:noncentered})
and $\fb\leq\fu$ 
(see \cite[Notes to Section 3]{vD:integer}). 

An ultrafilter $\mathcal{U}$ is called a \emph{p-point} 
if for every set $\{X_n\st n<\omega\}\subseteq\mathcal{U}$ 
there is an $X\in\mathcal{U}$ 
such that $X\aincl X_n$ for all $n<\omega$. 
The existence of a p-point 
is not provable in ZFC 
(\cite[Theorem 4.4.7]{BaJ:set}, see also \cite[Section VI.4]{Sh:improper}). 
For an uncountable regular cardinal $\kappa$, 
we say an ultrafilter $\mathcal{U}$ is a \emph{simple p${}_\kappa$-point} 
if it is generated by a decreasing sequence of length $\kappa$ 
with respect to $\aincl$. 
Clearly, for any $\kappa$, 
a simple p${}_\kappa$-point is a p-point. 

\begin{defn}
$\fpp$ is the smallest cardinal $\kappa$ 
for which a simple p${}_\kappa$-point exists 
(if such a $\kappa$ exists; otherwise $\fpp=\infty$). 
\end{defn}

Clearly we have $\fu\leq\fpp$, 
and $\fpp\leq\fc$ unless $\fpp=\infty$. 
Under Martin's axiom, 
there is a simple p${}_\fc$-point (\cite[Theorem 6.4]{FrZb:gaps}) 
and hence $\fpp=\fc$. 

Finally we define a cardinal invariant $\fl$, 
which was originally introduced in \cite{Kada:phd} 
(see also \cite[Section 5]{Kada:gamecd}). 
Let $\mathcal{S}=\prod_{n<\omega}[\omega]^{\leq n+1}$. 
We call each element of $\mathcal{S}$ a \emph{slalom}. 
%
%

\begin{defn}
$\fl$ is the smallest cardinal $\kappa$ satisfying the following: 
For every $h\in\omega^\omega$ 
there is 
a set $\Phi\subseteq\mathcal{S}$ of slaloms of size $\kappa$ such that, 
for every $f\in\prod_{n<\omega}h(n)$ 
there is a $\varphi\in\Phi$ such that for all but finitely many $n<\omega$ 
we have $f(n)\in\varphi(n)$. 
\end{defn}



\begin{prop}
\begin{enumerate}
\item\label{item:covmcovnl}
	\textup{(\cite[Lemma 7.2.3]{BaJ:set})}
	$\cov\M\leq\fl$ and $\cov\N\leq\fl$. 
\item\label{item:lcofn}
    \textup{(\cite[Theorem 2.3.9]{BaJ:set})}
	$\fl\leq\cof\N$. 
\end{enumerate}
\end{prop}


\begin{remark}
We say 
a proper forcing notion $\mathbb{P}$ 
has the \emph{Laver property} 
if for every $h\in\omega^\omega\cap\V$ 
and $g\in(\prod_{n<\omega}h(n))\cap\V^{\mathbb{P}}$ 
there is a $\varphi\in\mathcal{S}\cap\V$ 
such that 
for all $n<\omega$ we have $g(n)\in\varphi(n)$ 
(where $\V$ is a ground model). 
If $\V$ satisfies CH (the continuum hypothesis) 
and $\mathbb{P}$ has the Laver property, 
then $\fl=\aleph_1$ holds in $\V^{\mathbb{P}}$. 
\end{remark}

%
%


\section{Smirnov Compactifications}\label{sec:smirnov}


For a metric space $(X,d)$, 
$U^*_{d}(X)$ denotes the set of all bounded uniformly continuous functions 
from $(X,d)$ to $\R$.
$U^*_{d}(X)$ is a regular subring of $C^*(X)$. 
The \emph{Smirnov 
compactification $\smirnov{X}{d}$ of\/ $(X,d)$} 
is the unique compactification associated with the subring $U^*_{d}(X)$.
%
Note that, for disjoint closed subsets $A, B$ of $X$,
$\cptsep{A}{B}{u_d X}$ if and only if $d(A,B)>0$
(\cite[Theorem 2.5]{Wo:unifcpt}).

The following theorem tells us that 
we can approximate 
the Stone-\v{C}ech compactification of a metrizable space 
by the collection of all Smirnov compactifications. 
For a metrizable space $X$, 
let $\Metric(X)$ denote the set of all metrics on $X$ 
which are compatible with the topology on $X$. 

\begin{thm}\label{thm:smirnovapprox}
\textup{(\cite[Theorem 2.11]{Wo:unifcpt})}
For a non-compact metrizable space $X$, 
we have 
$\stonecech{X}\cptequiv\sup\{\smirnov{X}{d}\st d\in\Metric(X)\}$. 
\end{thm}

Now we set the following general question: 
\begin{quote}
\emph{How many metrics 
do we actually need to approximate the Stone-\v{C}ech compactification 
by Smirnov compactifications?}
\end{quote}
This question suggests the following definition of a cardinal function.

\begin{defn}
For a non-compact 
metrizable space $X$, 
$\fsa(X)$ is the smallest cardinality of a set $D\subseteq\Metric(X)$ 
which satisfies 
$\stonecech{X}\cptequiv\sup\{\smirnov{X}{d}\st d\in D\}$. 
\end{defn}

Let us begin with a familiar space: the half-open
interval $[0, \infty)$ equipped with the usual topology.

\begin{ex}\label{ex:sahalfopeninterval}
$\fsa([0,\infty))=\fd$. 
\end{ex}

\begin{proof}
First we prove that $\fsa([0,\infty))\leq\fd$. 
Take a set $\mathcal{F}\subseteq\omega^\omega$ 
which is of size $\fd$ 
and cofinal in $\omega^\omega$ with respect to $\domn$. 
We may assume that $f(n)\geq 1$ for all $f\in\mathcal{F}$ and $n<\omega$. 
For each $f\in\mathcal{F}$, 
define a metric $\rho_f$ on $[0,\infty)$ 
by letting 
\[
\rho_{f}(x,y)=\abs{\int_x^y f(\lfloor t\rfloor)dt}
\]
for $x,y\in[0,\infty)$ 
(where $\lfloor x\rfloor$ denotes the greatest integer not greater than $x$), 
and let $D=\{\rho_f\st f\in\mathcal{F}\}$. 
We will prove that $D$ satisfies 
$\stonecech{[0,\infty)}\cptequiv\sup\{\smirnov{[0,\infty)}{d}\st d\in D\}$. 

Take a pair of disjoint closed subsets $A,B$ of $[0,\infty)$. 
It suffices to show that $\cptsep{A}{B}{\smirnov{[0,\infty)}{\rho_f}}$ 
for some $f\in\mathcal{F}$. 
Without loss of generality, 
we may assume that 
$A$ and $B$ are of the form 
\[
A=\bigcup_{k<\omega}[c_{4k},c_{4k+1}]\text{ and }
B=\bigcup_{k<\omega}[c_{4k+2},c_{4k+3}] 
\]
for a non-decreasing sequence 
$\langle c_j\st j<\omega\rangle$ 
of non-negative real numbers 
such that $c_{2k+1}<c_{2k+2}$ for all $k<\omega$. 
Define $h_{A,B}\in\omega^\omega$ 
by letting 
\[
h_{A,B}(m)=
\max\left(\left\{\left\lceil\frac{1}{c_{2k+2}-c_{2k+1}}\right\rceil
	\st c_{2k+2}<m+1\right\}\cup\{1\}\right)
\]
for each $m<\omega$ 
(where $\lceil x\rceil$ denotes the smallest integer 
not smaller than $x$). 
Since $\mathcal{F}$ is cofinal in $\omega^\omega$ with respect to $\domn$, 
there is an $f\in\mathcal{F}$ which satisfies $h_{A,B}\domn f$. 
By the definition of $\rho_f$, 
for all but finitely many $k<\omega$ 
we have $\rho_f(c_{2k+1},c_{2k+2})\geq 1$, 
and hence $\rho_f(A,B)=\inf\{\rho_f(c_{2k+1},c_{2k+2})\st k<\omega\}>0$.
Thus we have $\cptsep{A}{B}{\smirnov{[0,\infty)}{\rho_f}}$.

We turn to the proof of $\fsa([0,\infty))\geq\fd$. 
Fix $\kappa<\fd$ and $D\subseteq\Metric([0,\infty))$ of size $\kappa$. 
We show that 
$\stonecech{[0,\infty)}\not\cptequiv
	\sup\{\smirnov{[0,\infty)}{d}\st d\in D\}$. 

For each $d\in D$, 
define $g_d\in\omega^\omega$ 
by letting 
\[
g_d(m)=\min\left\{k<\omega\st
d\left(m,m+\frac{1}{k}\right)<\frac{1}{m+1}\right\}
\]
for $m<\omega$. 
For each nonempty finite subset $F$ of $D$, 
let $g_F=\max\{g_f\st f\in F\}$ 
(where $\max$ is the pointwise maximum). 
Since $\size{[D]^{<\omega}}=\size{D}=\kappa<\fd$, 
there is an $f\in\omega^\omega$ 
which satisfies 
$f\not\domn g_F$ for every nonempty finite subsets $F$ of $D$. 
We may assume that $f(n)\geq 2$ for all $n<\omega$. 

Let $A=\omega$ and $B=\{m+\frac{1}{f(m)}\st m<\omega\}$. 
$A$, $B$ are disjoint closed subsets of $[0,\infty)$. 

For a nonempty finite subset $F$ of $D$, 
the set $I_F=\{m<\omega\st g_F(m)<f(m)\}$ 
is an infinite subset of $\omega$. 
Let $C=\cl\langle\bigcup\{U^*_d([0,\infty))\st d\in F\}\rangle$. 
$C$ is the closed subring of $C^*([0,\infty))$ 
associated with $\sup\{\smirnov{[0,\infty)}{d}\st d\in F\}$. 
By the definition of $g_F$, 
each $m\in I_F$ satisfies 
$d(m,m+\frac{1}{f(m)})<\frac{1}{m+1}$ for all $d\in F$.
If $\psi\in\bigcup\{U^*_d([0,\infty))\st d\in F\}$, 
then the sequence $\langle\psi(m)-\psi(m+\frac{1}{f(m)})\st m\in I_F\rangle$ 
converges to $0$. 
So for all $\varphi\in C$, 
$\langle\varphi(m)-\varphi(m+\frac{1}{f(m)})\st m\in I_F\rangle$ 
converges to $0$. 
This means that 
there are no $\varphi\in C$ and $a,b\in\R$ 
such that $a>b$, 
$\varphi(x)\geq a$ for all $x\in A$, 
and $\varphi(x)\leq b$ for all $x\in B$. 
By Lemma \ref{lem:charcptsep}, 
this means $\cptnonsep{A}{B}{\sup\{\smirnov{[0,\infty)}{d}\st d\in F\}}$. 
Since $F$ is an arbitrary nonempty finite subset of $D$ 
and by Lemma \ref{lem:cptsepcompact}, 
we have $\cptnonsep{A}{B}{\sup\{\smirnov{[0,\infty)}{d}\st d\in D\}}$, 
and hence 
$\stonecech{[0,\infty)}
	\not\cptequiv\sup\{\smirnov{[0,\infty)}{d}\st d\in D\}$. 
\end{proof}

It is obvious that, 
if $X$ is a discrete space and $d$ is the discrete metric on $X$ 
(that is, $d(x,y)=1$ whenever $x\neq y$), 
then 
$\stonecech{X}\cptequiv\smirnov{X}{d}$ 
and hence $\fsa(X)=1$.
%
In fact, $\fsa(X)=1$ holds if and only if the set of nonisolated points
of $X$ is compact (\cite[Corollary 3.5]{Wo:unifcpt}).

So it makes no sense to deal with $\fsa(\omega)$. 
Here we consider ``nontrivial'' ways to approximate $\beta\omega$ 
by Smirnov compactifications of $\omega$. 

For a metrizable space $X$, 
let $\Metric'(X)$ 
be the set of metrics $d\in\Metric(X)$ 
for which $\stonecech{X}\not\cptequiv\smirnov{X}{d}$. 

\begin{defn}
We say $D\subseteq\Metric'(\omega)$ satisfies 
the \emph{Smirnov finite intersection property} (Smirnov-FIP) 
if 
for every finite set $F\subseteq D$ 
we have 
$\stonecech{\omega}\not\cptequiv\sup\{\smirnov{\omega}{d}\st d\in F\}$. 
$\fsp$ is the smallest cardinality of a set $D\subseteq\Metric'(\omega)$ 
such that 
$D$ satisfies Smirnov-FIP 
and $\stonecech{\omega}\cptequiv\sup\{\smirnov{\omega}{d}\st d\in D\}$. 
\end{defn}

For a metrizable space $X$ and metrics $d_1,d_2\in\Metric(X)$, 
we write $d_1\preceq d_2$ 
if $U^*_{d_1}(X)\subseteq U^*_{d_2}(X)$ 
(or equivalently, $\smirnov{X}{d_1}\leq\smirnov{X}{d_2}$). 
Note that 
$d_1\preceq d_2$ 
if and only if 
the identity map on $X$ is uniformly continuous as a function 
from $(X,d_2)$ to $(X,d_1)$. 
If $d_1\leq d_2$ (that is, 
$d_1(x,y)\leq d_2(x,y)$ for all $x,y\in X$), 
then clearly $d_1\preceq d_2$. 

\begin{defn}
$\fsp'$ is the smallest cardinality of a set $D\subseteq\Metric'(\omega)$ 
such that 
$D$ is directed with respect to $\preceq$ 
(that is, 
for any $d_1,d_2\in D$ 
there is a $d\in D$ with $d_1\preceq d$ and $d_2\preceq d$)
and $\stonecech{\omega}\cptequiv\sup\{\smirnov{\omega}{d}\st d\in D\}$. 
\end{defn}

\begin{defn}
$\fst$ is the smallest cardinality of a set $D\subseteq\Metric'(\omega)$ 
such that 
$D$ is well-ordered by $\preceq$ 
and $\stonecech{\omega}\cptequiv\sup\{\smirnov{\omega}{d}\st d\in D\}$. 
\end{defn}

It is clear that $\fsp\leq\fsp'\leq\fst$.

\section{Analytic subgroups of $(2^\omega,+)$}\label{sec:analyticsubgroups}

In this section we give a lower bound for $\fsp$. 

For $f,g\in 2^\omega$, 
define $f+g\in 2^\omega$ by pointwise addition modulo 2. 
Then $(2^\omega,+)$ is an abelian group. 
The ideals $\M$ and $\N$ are invariant for ``translations'', 
that is, 
for any $A\subseteq 2^\omega$ and $f\in 2^\omega$, 
$A\in\M$ (resp. $\N$) 
if and only if 
$A+f=\{a+f\st a\in A\}\in\M$ (resp. $\N$).

\begin{lem}\label{lem:setofseparatedsets}
Suppose that 
$\alpha\omega$ is a compactification of $\omega$ 
and $A=2^\omega\cap C_{\alpha\omega}$. 
Then, 
\begin{enumerate}
\item\label{item:continuousextension} 
	For $B\subseteq\omega$, 
	$\cptsep{B}{\omega\ssm B}{\alpha\omega}$ 
    if and only if $\chi_B\in A$, where $\chi_B$ denotes the
    characteristic function of $B$.
\item\label{item:full} 
	$A=2^\omega$ 
	if and only if $\alpha\omega\cptequiv\stonecech{\omega}$. 
\item\label{item:subgroup} 
	$A$ is a subgroup of $(2^\omega,+)$. 
\item\label{item:tailset} 
	$A$ is closed under finite modifications, 
	that is, 
	if $f,g\in 2^\omega$ 
	and $f(n)=g(n)$ for all but finitely many $n<\omega$, 
	then $f\in A$ if and only if $g\in A$. 
\end{enumerate}
\end{lem}

\begin{proof}
(\ref{item:continuousextension}) is derived from Lemma \ref{lem:charcptsep}.
(\ref{item:full}), 
(\ref{item:subgroup}) 
and  
(\ref{item:tailset}) 
are easily seen. 
\end{proof}

Recall that a set in a Polish space $\mathcal{X}$ is called 
an \emph{analytic set} 
if it is the continuous image of a closed set in some Polish space, 
or equivalently, 
if it is the projection of a Borel set in $\mathcal{X}\times\mathcal{Y}$ 
for some Polish space $\mathcal{Y}$. 
Analytic sets are also called $\SIGMA^1_1$ sets. 
It is known that every analytic set has the Baire property 
and is Lebesgue measurable (\cite[Theorem 9.1.3]{BaJ:set}).

\begin{defn}
Let $\mathcal{A}$ 
be the collection of proper subgroups of $(2^\omega,+)$ 
which are analytic (in the Cantor space $2^\omega$) 
and closed under finite modifications. 
\end{defn}

\begin{lem}\label{lem:amcapn}
$\mathcal{A}\subseteq\M\cap\N$. 
\end{lem}

\begin{proof}
Follows from 
the zero-one law (see \cite[Theorems 8.47 and 17.1]{Ke:descriptive}) 
and the invariance of $\M$ and $\N$ for translations. 
\end{proof}


\begin{lem}\label{lem:unifcontiborel}
For $d\in\Metric(\omega)$, 
$U^*_d(\omega)$ is a $\PI^0_3$ set of the space $\R^\omega$. 
\end{lem}

\begin{proof}
$C^*(\omega)$ is a $\SIGMA^0_2$ subset of $\R^\omega$ 
because it is described 
by the following: 
\[
f\in C^*(\omega)
\iff
\exists m<\omega\,
\forall x<\omega\,
(f(x)\leq m).
\]
$U^*_d(\omega)$ is described by the following: 
\[
f\in U^*_d(\omega)
\iff
f\in C^*(\omega)
\,\land\,
\left[
\begin{gathered}
\forall\varepsilon\in\Q^+\,
\exists\delta\in\Q^+\,
\forall x,y<\omega\,\\
(d(x,y)\leq\delta
\to
\abs{f(x)-f(y)}\leq\varepsilon)
\end{gathered}
\right],
\]
where $\Q^+$ denotes the set of positive rational numbers. 
The latter formula of the right-hand side 
is a $\PI^0_3$ condition for $f\in\R^\omega$ 
and so is the whole right-hand side. 
Hence $U^*_d(\omega)$ is a $\PI^0_3$ set. 
\end{proof}

\begin{lem}\label{lem:finmetricanalytic}
Let $F=\{d_i\st i<n\}\subseteq\Metric(\omega)$ 
and 
$A_F
	=2^\omega\cap C_{\sup\{\smirnov{\omega}{d_i}\st i<n\}}$. 
Then $A_F$ is an analytic set in the Cantor space $2^\omega$. 
\end{lem}

\begin{proof}
Let $R=\langle\bigcup\{U^*_{d_i}(\omega)\st i<n\}\rangle$ 
and 
$C=\cl R$. 
Then $C=C_{\sup\{\smirnov{\omega}{d_i}\st i<n\}}$ 
and $A_F=2^\omega\cap C$. 

Suppose that $f\in A_F$. 
Then 
there is a sequence of functions in $R$ 
which uniformly converges to $f$. 
So there is a $g\in R$ such that, 
\begin{equation}\label{eq:unifconvergence}
\text{for all }x<\omega, 
\begin{cases}
g(x)\geq\frac{2}{3}	& \text{if }f(x)=1, \\
g(x)\leq\frac{1}{3}	& \text{if }f(x)=0.
\end{cases}
\tag{\textsf{U}${}_{f,g}$}
\end{equation}
On the other hand, 
suppose that there is a $g\in R$ satisfying \ref{eq:unifconvergence} 
for some $f\in 2^\omega$. 
Since $C=\cl R$ is a ring associated with a compactification of $\omega$
and by Lemma \ref{lem:charcptsep},
we have $f\in A_F$. 
Hence 
\[
A_F=\{	f\in 2^\omega\st
	\text{there is a }g\in R
	\text{ which satisfies \ref{eq:unifconvergence}}\}. 
\]
Since $R$ is a ring generated by 
the union of $U^*_{d_i}(\omega)$'s, 
every member $g$ of $R$ has the form 
$g=\sum_{k<m}\prod_{i<n}g_{k,i}$
for some $m<\omega$ and $g_{k,i}\in U^*_{d_i}(\omega)$ for each $k<m$, $i<n$. 
So $A_F$ is described as follows: 
$A_F=\bigcup\{A_F^m\st m<\omega\}$, and for each $m<\omega$, 
\begin{equation*}
A_F^m= \left\{
\begin{array}{c|c}
f\in 2^\omega &
\begin{gathered}
\text{there are }
g_{k,i}\in U^*_{d_i}(\omega)\text{ for }k<m\text{, }i<n\\
\text{such that }g=\textstyle\sum_{k<m}\prod_{i<n}g_{k,i}
\text{ satisfies \ref{eq:unifconvergence}} 
\end{gathered}
\end{array}
\right\}
\end{equation*}
Since the class of analytic sets is closed under countable unions, 
it suffices to prove that each $A_F^m$ is analytic. 
Fix $m<\omega$ 
and define a subset $B_F^m$ of the Polish space 
$2^\omega\times(\R^\omega)^{mn}$ 
by the following: 
\begin{equation*}
B_F^m= \left\{
\begin{array}{c|c}
\op{f}{\langle g_{k,i}\st k<m\text{, }i<n\rangle}&
\begin{gathered}
g_{k,i}\in U^*_{d_i}(\omega)\text{ for }k<m\text{, }i<n\\
\text{and }g=\textstyle\sum_{k<m}\prod_{i<n}g_{k,i}
\text{ satisfies \ref{eq:unifconvergence}} 
\end{gathered}
\end{array}
\right\}
\end{equation*}
Then $A_F^m$ is the projection of $B_F^m$ 
to the first coordinate. 
We check that $B_F^m$ is a Borel set. 
For each $k<m$ and $i<n$, 
the condition ``$g_{k,i}\in U^*_{d_i}(\omega)$'' 
is a $\PI^0_3$ condition by Lemma \ref{lem:unifcontiborel},
and the condition 
``$g=\sum_{k<m}\prod_{i<n}g_{k,i}$
	 satisfies \ref{eq:unifconvergence}'' 
is a closed ($\PI^0_1$) condition. 
So $B_F^m$ is a $\PI^0_3$ (and hence Borel) set of 
$2^\omega\times(\R^\omega)^{mn}$. 
\end{proof}

\begin{thm}\label{thm:covasp}
$\fsp\geq\cov{\mathcal{A}}$. 
\end{thm}

\begin{proof}
Suppose that 
$D\subseteq\Metric'(\omega)$ 
satisfies the Smirnov-FIP 
and $\sup\{\smirnov{\omega}{d}\st d\in D\}\cptequiv\stonecech{\omega}$. 
For a nonempty finite set $F\subseteq D$, 
let 
$A_F
	=2^\omega\cap C_{\sup\{\smirnov{\omega}{d}\st d\in F\}}$. 
By the Smirnov-FIP of $D$ and Lemma \ref{lem:setofseparatedsets}, 
$A_F$ is a proper subgroup of $(2^\omega,+)$ 
and closed under finite modifications. 
Also, 
by Lemma \ref{lem:finmetricanalytic}, 
$A_F$ is an analytic set in the Cantor space $2^\omega$. 
Hence we have $A_F\in\mathcal{A}$ for each $F$. 

Since $\sup\{\smirnov{\omega}{d}\st d\in D\}\cptequiv\stonecech{\omega}$ 
and by Lemma \ref{lem:cptsepcompact},
for each $B\subseteq\omega$ 
there is a finite set $F\subseteq D$ 
for which 
$\cptsep{B}{\omega\ssm B}{\sup\{\smirnov{\omega}{d}\st d\in F\}}$, 
and by Lemma \ref{lem:setofseparatedsets} 
and the definition of $A_F$, 
$\chi_B$ belongs to $A_F$. 
This means that 
the set 
$\{A_F\st F\in[D]^{<\omega}\ssm\{\emptyset\}\}$ 
covers $2^\omega$. 
Since each $A_F$ belongs to $\mathcal{A}$, 
we have 
$\size{[D]^{<\omega}}=\size{D}\geq\cov{\mathcal{A}}$. 
\end{proof}

\begin{cor}
$\fsp\geq\cov\M$ and $\fsp\geq\cov\N$. 
\end{cor}

\begin{proof}
Follows from Lemma \ref{lem:amcapn} and Theorem \ref{thm:covasp}. 
\end{proof}

\section{Ultrafilters and slaloms}\label{sec:uf}

In this section we investigate upper bounds for $\fsp'$ and $\fst$ 
using infinitary combinatorics on $\omega$. 
First we show that an ultrafilter gives 
an upper bound for $\fsp'$. 

\begin{thm}\label{thm:spu}
$\fsp'\leq\fu$. 
\end{thm}

\begin{proof}
Let $\mathcal{U}$ be a subset of $\Pow{\omega}$ of size $\fu$ 
which generates an ultrafilter. 
For each $X\in\mathcal{U}$, 
define a metric $d_X$ on $\omega$ by the following: 
\[
d_X(x,y)=
\begin{cases}
0	&	\text{if $x=y$},	\\
\abs{2^{-x}-2^{-y}}	&	\text{if $x\neq y$ and $x,y\in X$},	\\
1	&	\text{otherwise}.
\end{cases}
\]
Note that, 
for $B\subseteq\omega$,
$d_X(B, \omega\ssm B)>0$ (or equivalently
$\cptsep{B}{\omega\ssm B}{\smirnov{\omega}{d_X}}$) holds
if and only if $X\aincl B$ or $X\aincl\omega\ssm B$ holds.

Let $D=\{d_X\st X\in\mathcal{U}\}$. 
Clearly we have $D\subseteq\Metric'(\omega)$, 
and for $X,Y\in\mathcal{U}$, 
$X\subseteq Y$ implies $d_Y\leq d_X$ 
and $X\aincl Y$ implies $d_Y\preceq d_X$. 
Since $\mathcal{U}$ is directed with respect to $\supseteq$, 
$D$ is directed with respect to $\preceq$. 


Since $\mathcal{U}$ generates an ultrafilter, 
for each $B\subseteq\omega$ we can find $X\in\mathcal{U}$ 
so that 
$X\aincl B$ or $X\aincl\omega\ssm B$. 
This means that 
$\cptsep{B}{\omega\ssm B}{\sup\{\smirnov{\omega}{d}\st d\in D\}}$ 
for all $B\subseteq\omega$, 
and hence 
$\stonecech{\omega}\cptequiv\sup\{\smirnov{\omega}{d}\st d\in D\}$. 
\end{proof}



\begin{cor}\label{cor:stpp}
$\fst\leq\fpp$. 
\end{cor}

\begin{proof}
Modify the proof of Theorem \ref{thm:spu}. 
\end{proof}

Another upper bound for $\fsp'$ 
is given by slaloms. 

\begin{thm}\label{thm:spl}
$\fsp'\leq\fl$. 
\end{thm}

\begin{proof}
Fix a partition 
$\mathcal{I}=\{I_n\st n<\omega\}$ of $\omega$ 
such that 
$\size{I_n}=2^{(n+1)^2}$ for all $n<\omega$. 
We say 
a partition $\mathcal{P}$ of $\omega$ 
is a \emph{moderate refinement of\/ $\mathcal{I}$\/} 
if $\mathcal{P}$ refines $\mathcal{I}$ 
and 
each $I_n$ is partitioned into at most $2^{n+1}$ pieces in $\mathcal{P}$. 
Let $\mathcal{R}$ denote 
the collection of all moderate refinements of $\mathcal{I}$. 

For each finite 
subset $\mathcal{F}$ of $\mathcal{R}$, 
define a metric $d_{\mathcal{F}}$ on $\omega$ as follows: 
\[
d_{\mathcal{F}}(x,y)=
\begin{cases}
0	&	\text{if $x=y$,}	\\
\frac{1}{n+1}	
	&	\text{if $x\neq y$, $x,y\in I_n$ and 
		$\forall\mathcal{P}\in\mathcal{F}\,
		\exists B\in\mathcal{P}\,(x,y\in B)$,}
	\\
1	&	\text{otherwise.}
\end{cases}
\]
It is easy to see that 
$d_{\mathcal{F}}\in\Metric(\omega)$. 
Moreover, 
for $\mathcal{F},\mathcal{G}\in[\mathcal{R}]^{<\omega}$, 
$\mathcal{F}\subseteq\mathcal{G}$ implies 
$d_{\mathcal{F}}\leq d_{\mathcal{G}}$, 
and hence 
$D=\{d_{\mathcal{F}}\st\mathcal{F}\in[\mathcal{R}]^{<\omega}\}$
is directed with respect to $\preceq$. 

\begin{claim}
$D\subseteq\Metric'(\omega)$. 
\end{claim}

\begin{proof}
It suffices to show that 
for every $\mathcal{F}\in[\mathcal{R}]^{<\omega}$ 
and $\varepsilon>0$ 
there are distinct $x,y<\omega$ 
such that $d_{\mathcal{F}}(x,y)<\varepsilon$. 
Let $k=\size{\mathcal{F}}$, 
and let $\mathcal{P}_{\mathcal{F}}$ 
be the coarsest common refinement of the partitions in $\mathcal{F}$. 
Then, 
for each $n<\omega$, 
$I_n$ is partitioned 
into at most $2^{(n+1)k}$ pieces in $\mathcal{P}_{\mathcal{F}}$. 
Hence, 
for $n\geq k$, 
there are distinct $x,y\in I_n$ 
contained in the same piece in $\mathcal{P}_{\mathcal{F}}$, 
which means that 
$d_{\mathcal{F}}(x,y)=\frac{1}{n+1}$. 
\end{proof}

By the definition of $\fl$, 
we choose a set $\Phi\subseteq\prod_{n<\omega}[2^{I_n}]^{\leq n+1}$ 
of size $\fl$ 
which satisfies the following: 
for every $f\in 2^\omega$ 
there is a $\varphi\in\Phi$ 
such that 
for all $n<\omega$ we have $f\restrictedto I_n\in\varphi(n)$. 
For each $\varphi\in\Phi$, 
define a moderate refinement $\mathcal{P}_\varphi$ of $\mathcal{I}$ 
by the following: 
$x,y\in I_n$ are contained in the same piece in $\mathcal{P}_\varphi$ 
if and only if 
$s(x)=s(y)$ for all $s\in\varphi(n)$. 
Let $\mathcal{R}_\Phi=\{\mathcal{P}_\varphi\st\varphi\in\Phi\}$ 
and $D_\Phi=\{d_{\mathcal{F}}\st\mathcal{F}\in[\mathcal{R}_\Phi]^{<\omega}\}$. 
Then $\size{D_\Phi}=\size{[\Phi]^{<\omega}}=\size{\Phi}=\fl$, 
$D_\Phi\subseteq\Metric'(\omega)$ 
and $D_\Phi$ is directed with respect to $\preceq$. 

\begin{claim}
$\sup\{\smirnov{\omega}{d}\st d\in D_\Phi\}\cptequiv\stonecech{\omega}$. 
\end{claim}

\begin{proof}
It suffices to show that 
every $f\in 2^\omega$ belongs to $U^*_d(\omega)$ 
for some $d\in D_\Phi$. 
Fix $f\in 2^\omega$. 
By the choice of $\Phi$, 
there is a $\varphi\in\Phi$ such that 
for all $n<\omega$ we have $f\restrictedto I_n\in\varphi(n)$. 
Let $d=d_{\{\mathcal{P}_\varphi\}}\in D_\Phi$. 
Then for every $x,y<\omega$, 
$d(x,y)=1$ whenever $f(x)\neq f(y)$ holds, 
and hence $f\in U^*_d(\omega)$. 
\end{proof}

This concludes the proof. 
\end{proof}


\section{Consistency results}\label{sec:con}

In this section 
we observe various consistency results 
concerning $\fsp$, $\fsp'$, $\fst$ 
and other cardinal invariants of the reals. 

Suppose that $\kappa\geq\aleph_2$ and $\cf(\kappa)\geq\aleph_1$. 
When we add $\kappa$-many Cohen reals 
over a model satisfying CH (the continuum hypothesis), 
we get a model which satisfies 
$\aleph_1=\non\M<\cov\M=\fc=\kappa$. 
On the other hand, 
if we add $\kappa$-many random reals over a model for CH, 
$\aleph_1=\non\N=\fd<\cov\N=\fc=\kappa$ holds in the forcing model 
(see \cite[Chapter 3]{BaJ:set} or \cite{Ku:reals}). 
So we have the following consistency results.

\begin{thm}\label{thm:condsp}
Each of the following statements is consistent with 
$\operatorname{ZFC}$\textup{:}
\begin{enumerate}
\item $\aleph_1=\non\M<\fsp$. 
\item $\aleph_1=\non\N=\fd<\fsp$. 
\end{enumerate}
\end{thm}

Hence none of $\non\M$, $\non\N$ and $\fd$ 
can be an upper bound for $\fsp$. 

A model of ZFC satisfying $\aleph_1=\fl<\fb=\non\N=\fc=\aleph_2$ 
is obtained by a countable support iteration of Mathias forcing 
of length $\omega_2$ over a model 
satisfying CH (\cite[Subsection 7.4.A]{BaJ:set}). 

\begin{thm}\label{thm:conspb}
$\aleph_1=\fsp'<\fb=\non\N=\aleph_2$ is consistent 
with $\operatorname{ZFC}$. 
\end{thm}

We already know that $\cov\M$ and $\cov\N$ are lower bounds for $\fsp$. 
The above theorem means that 
no other (nontrivial) lower bound for $\fsp$ is found in Cicho\'n's diagram 
(cf.\ Figure \ref{fig:cd}). 



We turn to the consistency results involving $\fst$. 

CH implies the existence of a simple p${}_{\aleph_1}$-point. 
A countable support iteration of 
the infinitely equal forcing (\cite[Subsection 7.4.C]{BaJ:set}) 
preserves p-points 
(that is, 
a p-point in the ground model still generates an ultrafilter 
in the forcing model). 
So $\fpp=\aleph_1$ holds in the forcing model 
obtained by 
a countable support iteration of the infinitely equal forcing 
of length $\omega_2$ 
over a model for CH. 
On the other hand, 
it is easily checked that 
$\fl=\fc=\aleph_2$ holds in the same model. 

\begin{thm}
$\aleph_1=\fsp=\fsp'=\fst<\fl=\fc=\aleph_2$ is consistent 
with $\operatorname{ZFC}$. 
\end{thm}

Martin's axiom 
implies 
$\cov\M=\cov\N=\fc$ 
and the existence of a simple p${}_\fc$-point, 
and hence $\fsp=\fsp'=\fst=\fpp=\fc$. 

\begin{thm}
$\aleph_1<\fsp=\fsp'=\fst=\fc$ is consistent 
with $\operatorname{ZFC}$. 
\end{thm}

Now we are going to prove that 
$\fst=\infty$ is consistent with ZFC, 
that is, 
possibly 
we cannot approximate $\stonecech{\omega}$ 
by any $\leq$-increasing chain of Smirnov compactifications of $\omega$. 

We use the following theorem, which is due to Kunen. 

\begin{thm}\label{thm:isomorphismofnames}
\textup{(\cite{Ku:inacc})}
The following holds 
in the forcing model obtained by adding $\aleph_2$-many Cohen reals 
over a model for $\operatorname{CH}$: 
Let $\mathcal{X}$ be a Polish space 
and $A\subseteq\mathcal{X}\times\mathcal{X}$ a Borel set. 
Then there is no sequence $\{r_\alpha\st\alpha<\omega_2\}$ 
in $\mathcal{X}$ 
which satisfies 
\[
\alpha\leq\beta<\omega_2 
\text{ if and only if } \op{r_\alpha}{r_\beta}\in A.	
\]
\end{thm}

An easy consequence of Kunen's theorem is that, 
in the forcing model obtained by adding $\aleph_2$-many Cohen reals 
over a model for CH, 
there is no strictly $\aincl$-increasing (or decreasing) 
sequence of length $\omega_2$ in $\Pow{\omega}$. 
On the other hand, 
$\cov\M=\fu=\fc=\aleph_2$ holds in this model, 
and hence $\fpp=\infty$. 

We may regard a metric on $\omega$ 
as an element of the Polish space $\R^{\omega\times\omega}$. 
We define $A\subseteq(\R^{\omega\times\omega})^2$ 
by the following: 
\[
\op{f_1}{f_2}\in A 
	\iff
	\left[
	\begin{gathered}
	\forall\varepsilon\in\Q^+\,
	\exists\delta\in\Q^+\,
	\forall\op{x}{y}\in\omega\times\omega\,	\\
	(f_2(x,y)\geq\delta
	\,\lor\, f_1(x,y)\leq\varepsilon)
	\end{gathered}
	\right].
\]
Then $A$ is a $\PI^0_3$ 
subset of $(\R^{\omega\times\omega})^2$ and, 
for $d_1,d_2\in\Metric(\omega)$, 
$\op{d_1}{d_2}\in A$ 
if and only if 
the identity map on $\omega$ is uniformly continuous as a function 
from $(\omega,d_2)$ to $(\omega,d_1)$, 
that is, $d_1\preceq d_2$. 
So Theorem \ref{thm:isomorphismofnames} leads the following lemma. 

\begin{lem}
In the forcing model obtained by adding $\aleph_2$-many Cohen reals 
over a model for $\operatorname{CH}$, 
there is no sequence $\{d_\alpha\st\alpha<\omega_2\}$ in $\Metric'(\omega)$ 
which is strictly $\preceq$-increasing, 
that is, 
$d_\alpha\preceq d_\beta$ if and only if $\alpha\leq\beta$. 
\end{lem}

On the other hand, 
$\cov\M=\fu=\fc=\aleph_2$ and hence $\fsp=\aleph_2$ 
holds in this model. 
So we have the following consistency result. 

\begin{thm}
$\fst=\infty$ is consistent with $\operatorname{ZFC}$. 
\end{thm}

\section{Higson Compactifications}\label{sec:higson}

For a metrizable space $X$, 
a metric $d$ on $X$ is called \emph{proper} 
if each $d$-bounded set has compact closure. 
A \emph{proper metric space} means 
a metric space whose metric is proper. 
For a metrizable space $X$, 
let $\PM(X)$ be the set of all proper metrics 
compatible with the topology of $X$. 
It is known that, 
for every locally compact separable metrizable space $X$, 
we have $\PM(X)\neq\emptyset$ (\cite[Lemma 3.1]{KT:approx}). 

Let $(X,d)$ be a proper metric space and $(Y,\rho)$ a metric space. 
We say 
a function 
$f$ from $X$ to $Y$ is \emph{slowly oscillating} 
if it satisfies the following condition: 
\begin{equation}
\begin{gathered}
\forall r>0\,
\forall\varepsilon>0\, 
\exists K\subseteq X\text{ a compact set }\, \\
\forall x\in X\ssm K\,(\diam_\rho(f''\Ball_d(x,r))<\varepsilon). 
\end{gathered}
\tag*{$(*)_d$}
\end{equation}
For a proper metric space $(X,d)$, 
let $C^*_d(X)$ be 
the set of all bounded continuous slowly oscillating functions 
from $(X,d)$ to $\R$.
$C^*_d(X)$ is a regular subring of $C^*(X)$.
The \emph{Higson compactification $\higson{X}{d}$ of $(X,d)$} 
is the unique compactification associated with the subring $C^*_d(X)$.
%
%
Note that, for disjoint closed subsets $A,B$ of $X$,
    $\cptsep{A}{B}{\higson{X}{d}}$
    if and only if
    for any $R>0$ there is a compact subset $K_R$ of $X$ 
    such that 
    $\sum_{i<n}d(x,E_i)>R$ holds for all $x\in X\ssm K_R$.
    (See also \ \cite[Proposition 2.3]{DKU:higsoncorona}.)

The following 
is a corresponding proposition of Theorem \ref{thm:smirnovapprox} 
for Higson compactifications. 

\begin{thm}\label{thm:higsonapprox}
\textup{(\cite[Theorem 3.2]{KT:approx})}
For a non-compact locally compact separable metrizable space $X$, 
we have 
$\stonecech{X}\cptequiv\sup\{\higson{X}{d}\st d\in\PM(X)\}$. 
\end{thm}

So we consider the following cardinal function. 

\begin{defn}
For a non-compact locally compact separable metrizable space $X$, 
$\fha(X)$ is the smallest cardinality of a set $D\subseteq\PM(X)$ 
which satisfies 
$\stonecech{X}\cptequiv\sup\{\higson{X}{d}\st d\in D\}$. 
\end{defn}

\begin{lem}\label{lem:saha}
For a non-compact locally compact separable metric space $X$, 
we have $\fsa(X)\leq\fha(X)$. 
\end{lem}

\begin{proof}
It is clear that $\PM(X)\subseteq\Metric(X)$. 
For each $d\in\PM(X)$ we have $\higson{X}{d}\leq\smirnov{X}{d}$, 
and so 
if $D\subseteq\PM(X)$ satisfies 
$\stonecech{X}\cptequiv\sup\{\higson{X}{d}\st d\in D\}$, 
it also satisfies 
$\stonecech{X}\cptequiv\sup\{\smirnov{X}{d}\st d\in D\}$. 
\end{proof}

\begin{ex}\label{ex:hahalfopeninterval}
$\fha([0,\infty))=\fd$. 
\end{ex}

\begin{proof}

$\fha([0,\infty))\geq\fd$ 
follows from Lemma \ref{lem:saha} and Example \ref{ex:sahalfopeninterval}. 
The other inequality 
is proved 
in the same way as the proof of Example \ref{ex:sahalfopeninterval}, 
using 
\[
h'_{A,B}(m)
=\max\left(
	\left\{
		k\cdot
		\left\lceil
			\frac{1}{c_{2k+2}-c_{2k+1}}
		\right\rceil
		\st c_{2k+2}<m+1
	\right\}
	\cup\{1\}
\right)
\]
instead of $h_{A,B}$. 
(See also \cite[Section 2]{IT:halfopen}.)
\end{proof}

%
For a non-compact locally compact separable metrizable space $X$,
$\fha(X)=1$ holds if and only if the set of nonisolated points of $X$ is
compact (\cite[Proposition 2.6]{KT:approx}).

For a metrizable space $X$, 
let $\PM'(X)$ 
be the set of proper metrics $d\in\PM(X)$ 
for which $\stonecech{X}\not\cptequiv\higson{X}{d}$. 

\begin{defn}
We say $D\subseteq\PM'(\omega)$ satisfies 
the \emph{Higson finite intersection property} (Higson-FIP) 
if 
for every finite set $F\subseteq D$ 
we have 
$\stonecech{\omega}\not\cptequiv\sup\{\higson{\omega}{d}\st d\in F\}$. 
$\fhp$ is the smallest cardinality of a set $D\subseteq\PM'(\omega)$ 
such that 
$D$ satisfies Higson-FIP 
and $\stonecech{\omega}\cptequiv\sup\{\higson{\omega}{d}\st d\in D\}$. 
\end{defn}

For a space $X$ and metrics $d_1,d_2\in\Metric(X)$, 
we write $d_1\sqsubseteq d_2$ 
if $C^*_{d_1}(X)\subseteq C^*_{d_2}(X)$ 
(or equivalently, $\higson{X}{d_1}\leq\higson{X}{d_2}$). 

\begin{defn}
$\fhp'$ is the smallest cardinality of a set $D\subseteq\PM'(\omega)$ 
such that 
$D$ is directed with respect to $\sqsubseteq$ 
and $\stonecech{\omega}\cptequiv\sup\{\higson{\omega}{d}\st d\in D\}$. 
\end{defn}

\begin{defn}
$\fht$ is the smallest cardinality of a set $D\subseteq\PM'(\omega)$ 
such that 
$D$ is well-ordered by $\sqsubseteq$ 
and $\stonecech{\omega}\cptequiv\sup\{\higson{\omega}{d}\st d\in D\}$. 
\end{defn}

It is clear that $\fhp\leq\fhp'\leq\fht$. 
We can prove $\cov{\mathcal{A}}\leq\fhp\leq\fhp'\leq\fl$ 
by arguments similar to the ones in 
Sections \ref{sec:analyticsubgroups} and \ref{sec:uf}. 

\begin{lem}\label{lem:slowoscborel}
For a proper metric space $(\omega,d)$, 
$C^*_d(\omega)$ is a $\PI^0_3$ subset of $\R^\omega$. 
\end{lem}

\begin{proof}
The condition $(*)_d$ is 
equivalent to the following: 
\begin{align*}
&	\forall R,\varepsilon\in\Q^+\,
\exists k<\omega\,
\forall x,y<\omega\,	\\
&	(x\geq k\,\land\,y\geq k\,\land\,d(x,y)\leq R
\to
\abs{f(x)-f(y)}\leq\varepsilon).
\end{align*}
This formula is a $\PI^0_3$ condition for $f\in\R^\omega$. 
Since $C^*(\omega)$ is a $\SIGMA^0_2$ subset of $\R^\omega$, 
$C^*_d(\omega)$ is a $\PI^0_3$ subset of $\R^\omega$. 
\end{proof}

\begin{thm}\label{thm:covahp}
$\fhp\geq\cov{\mathcal{A}}$. 
\end{thm}

\begin{proof}
Proved in the same way 
as the proof of Theorem \ref{thm:covasp}, 
using Lemma \ref{lem:slowoscborel}
instead of Lemma \ref{lem:unifcontiborel}. 
\end{proof}

\begin{cor}
$\fhp\geq\cov\M$ and $\fhp\geq\cov\N$. 
\end{cor}

\begin{thm}\label{thm:hpl}
$\fhp'\leq\fl$. 
\end{thm}

\begin{proof}
Proved in the same way as the proof of Theorem \ref{thm:spl}, 
using the following metric for $\mathcal{F}\in[\mathcal{R}]^{<\omega}$ 
instead of $d_{\mathcal{F}}$. 
\[
\rho_{\mathcal{F}}(x,y)=
\begin{cases}
0	&	\text{if $x=y$,}	\\
1	&	
	\text{if $x\neq y$ and 
		$\forall\mathcal{P}\in\mathcal{F}\,
		\exists B\in\mathcal{P}\,(x,y\in B)$,}	\\
1+\max\{m,n\}	&
	\text{otherwise, and if $x\in I_m$, $y\in I_n$. }
\end{cases}
\] 
\end{proof}

We have the following consistency results 
concerning $\fhp$, $\fhp'$ and cardinals from Cicho\'n's diagram. 

\begin{thm}\label{thm:condhp}
\textup{(cf.\ Theorem \ref{thm:condsp})}
Each of the following statements is consistent 
with $\operatorname{ZFC}$\textup{:}
\begin{enumerate}
\item $\aleph_1=\non\M<\fhp$. 
\item $\aleph_1=\non\N=\fd<\fhp$. 
\end{enumerate}
\end{thm}

\begin{thm}\label{thm:conhpb}
\textup{(cf.\ Theorem \ref{thm:conspb})}
$\aleph_1=\fhp'<\fb=\non\N=\aleph_2$ is consistent 
with $\operatorname{ZFC}$. 
\end{thm}

Now we prove that $\fht=\infty$ is consistent with ZFC. 
We need the following characterization of the $\sqsubseteq$-relation, 
which appears in \cite{Iwa:manuscript,Tomo:higsonproduct}. 

\begin{lem}\label{lem:charhigsonorder}
For proper metric space $(X,d_1)$ and $(X,d_2)$, 
the following are equivalent\/\textup{:} 
\begin{enumerate}
\item $d_1\sqsubseteq d_2$. 
\item For all $r>0$ 
	there is a $\varepsilon>0$ such that, 
	for all $x,y\in X$ 
	if $d_2(x,y)<r$ 
	then $d_1(x,y)<\varepsilon$. 
\end{enumerate}
\end{lem}

\begin{lem}
In the forcing model obtained by adding $\aleph_2$-many Cohen reals 
over a model for $\operatorname{CH}$, 
there is no sequence $\{d_\alpha\st\alpha<\omega_2\}$ in $\PM'(\omega)$ 
which is strictly $\sqsubseteq$-increasing, 
that is, 
$d_\alpha\sqsubseteq d_\beta$ if and only if $\alpha\leq\beta$. 
\end{lem}

\begin{proof}
Define $A\subseteq(\R^{\omega\times\omega})^2$ 
by 
\[
\op{f_1}{f_2}\in A 
	\iff
	\left[
	\begin{gathered}
	\forall r\in\Q^+\,
	\exists\varepsilon\in\Q^+\,
	\forall\op{x}{y}\in\omega\times\omega\,	\\
	(f_2(x,y)\geq r
	\,\lor\, f_1(x,y)\leq\varepsilon)
	\end{gathered}
	\right]
\]
and apply Theorem \ref{thm:isomorphismofnames}. 
\end{proof}

In the same model, 
we have $\cov\M=\fhp=\fc=\aleph_2$. 

\begin{thm}
$\fht=\infty$ is consistent with $\operatorname{ZFC}$. 
\end{thm}

\section{Questions}\label{sec:q}

\begin{q}
$\fhp'\leq\fu$? \quad
$\fht\leq\fpp$?
\end{q}

\begin{q}
$\fsp=\fsp'$? \quad 
$\fhp=\fhp'$?
\end{q}

\begin{q}
\begin{enumerate}
\item Is $\cof\M<\fsp$ consistent with $\operatorname{ZFC}$?
\item Is $\cof\M<\fhp$ consistent with $\operatorname{ZFC}$?
\end{enumerate}
\end{q}

\begin{q}
\begin{enumerate}
\item Is $\fsp'<\fst\leq\fc$ consistent with $\operatorname{ZFC}$?
\item Is $\fhp'<\fht\leq\fc$ consistent with $\operatorname{ZFC}$?
\end{enumerate}
\end{q}

\begin{q}
Is it consistent with $\operatorname{ZFC}$
that $\fst\leq\fc$ and $\fpp=\infty$?
\end{q}

\subsection*{Acknowledgement}

The authors would like to thank
J\"org Brendle, Saka\'e Fuchino, and Hiroshi Fujita
for helpful suggestion and discussion 
during this work.


\vspace*{12pt}
\indent
Masaru Kada,
    Information Processing Center, 
	Kitami Institute of Technology. 
	Kitami 090--8507 JAPAN.\par
\indent
e-mail: kada@math.cs.kitami-it.ac.jp\par
\indent
Kazuo Tomoyasu,
	General Education, Miyakonojo National College of Technology,
	Miyakonojo-shi, Miyazaki 885-8567 JAPAN.\par
\indent
e-mail:tomoyasu@cc.miyakonojo-nct.ac.jp\par
\indent
Yasuo Yoshinobu,
	Graduate School of Information Science,	
	Nagoya University. 
	Nagoya 464--8601 JAPAN.\par
\indent
e-mail:yosinobu@math.nagoya-u.ac.jp

\end{document}